\documentclass[draft,11pt]{article}
\usepackage{amssymb,amsmath}
\usepackage{mathrsfs}
\usepackage[english]{babel}
\newtheorem{theo}{Theorem}[section]
\newtheorem{lema}[theo]{Lemma}
\newtheorem{nota}[theo]{Remark}
\newtheorem{definicion}[theo]{Definition}
\newtheorem{corolary}[theo]{Corollary}
\newenvironment{proof}{{\bf Proof:\ }}{\hfill
$\Box$}

\newcommand{\bq}{\begin{equation}}
\newcommand{\eq}{\end{equation}}
\newcommand{\ba}{\begin{array}}
\newcommand{\ea}{\end{array}}

\newcommand{\half}{\frac{1}{2}}

\newcommand{\dst}{\displaystyle}

\newcommand\be{\begin{enumerate}}
\newcommand\ee{\end{enumerate}}
\newcommand\bi{\begin{itemize}}
\newcommand\ei{\end{itemize}}

\newcommand\XX[1]{\mathbb{#1}}

\newcommand\bd{\begin{definicion}{\bf }}
\newcommand\ed{\end{definicion}}
\newcommand\bl{\begin{lema}{\bf }}
\newcommand\el{\end{lema}}
\newcommand\bp{\begin{prop}{\bf }}
\newcommand\ep{\end{prop}}
\newcommand\bt{\begin{theo}{\bf }}
\newcommand\et{\end{theo}}
\newcommand\bdm{\begin{proof}}
\newcommand\edm{\end{proof}}
\newcommand\bn{\begin{nota}{\bf }}
\newcommand\en{\end{nota}}
\newcommand\bc{\begin{corolary}{\bf }}
\newcommand\ec{\end{corolary}}

\newfont{\got}{eufm10 scaled \magstep1}

\title{On the elementary symmetric functions of a sum of matrices}
\date{\today}
\author{R. S. Costas-Santos\thanks{
EMAIL: rscosa@gmail.com, URL: www.rscosa.com \newline
Department of Mathematics, University of California, Santa
Barbara, CA 93106, US}}

\begin{document}
\maketitle

\begin{abstract}
Often in mathematics it is useful to summarize a
multivariate phenomenon with a single number.
In fact, the determinant -- which is denoted
by $\det$ -- is one of the simplest cases and
many of its properties are very well-known.
For instance, the determinant is a multiplicative
function, i.e.  $\det (A\, B)=\det A\cdot \det B$,
$A, B \in M_n$, and it is a multilinear
function, but it  is not, in general, an additive
function, i.e.  $\det (A+B)\ne \det A+\det B$.
\\
Another interesting scalar function in the Matrix
Analysis is the characteristic polynomial.
In fact, given a square matrix $A$, the coefficients
of  its characteristic polynomial $\chi_A(t):=\det(t I-A)$
are, up to a sign, the elementary symmetric functions
associated with the eigenvalues of $A$.
\\
In the present paper we present new expressions
related to the elementary symmetric functions of sum of
matrices.

The main motivation of this manuscript
is try to find new properties to probe
the following conjecture.

{\bf Bessis-Moussa-Villani conjecture:}
\cite{BeMoVi75,LiSe04}\\
{\it The polynomial $p(t):=Tr((A+tB)^m)\in \XX R[t]$,
has only nonnegative coefficients whenever $A, B\in M_r$
are positive semidefinite matrices.}

Moreover, some numerical evidences and the Newton-Girard
formulas suggested to us to consider a more general
conjecture that will be considered in a further manuscript.

{\bf Positivity Conjecture :} \\
{\it The polynomial $S_k((A+tB)^m)\in \XX R[t]$, has only
nonnegative coefficients whenever $A, B\in M_r$ are positive
semidefinite matrices for every $k=0,1, \dots,r.$}

It is clear that the BMV conjecture is a particular case of the
positivity conjecture for $k=1$ since $S_1\equiv Tr$.
\end{abstract}

\noindent{\em key words and phrases:} Elementary symmetric function,
Hermitian matrix, determinant.

\section{Introduction}
Denote by $M_{m,n}$ the set of $m\times n$ matrices
over an arbitrary field $\XX F$ and by $M_n$ the set
$M_{n,n}$.
Determinants are mathematical objects that are very
useful in the matrix analysis.

In fact, the determinant of a matrix $A\in M_n$,
can be presented in two important, apparently different,
but equivalent ways.

The first one is the {\sc Laplace expansion}:
\\
If $A=[a_{i,j}]$ and assuming that the determinant
is defined over $M_{n-1}$, then
\bq \label{ex5}
\det(A)=\sum_{j=1}^n(-1)^{i+j}a_{i,j}\det(A_{i,j})=
\sum_{i=1}^n(-1)^{i+j}a_{i,j}\det(A_{i,j}),
\eq
where $A_{i,j}\in M_{n-1}$ denotes the
submatrix  of $A$ resulting from the deletion
of row $i$ and  column $j$.
\\
The second way is the {\sc alternating sum}:
\bq \label{ex6}
\det(A)=\sum_{\sigma\in P_n}{\rm sgn}(\sigma)\,
a_{1,\sigma(1)} a_{2,\sigma(2)}\cdots a_{n,\sigma(n)},
\eq
where $P_n$ is the set of all permutation of
$\{1,2,\dots,n\}$, and ${\rm sgn}(\sigma)$ denotes the
sign of the permutation $\sigma$.

\bn
Notice that with these definitions it is clear
that the determinant is a multilinear function.
\en
In the present paper we present a closed expression for
$$
\det(A_1+A_2+\cdots+A_N),\qquad A_1,A_2,\dots,A_N\in M_n,
$$
where $N\ge n+1$, in terms of the sum of another
determinants involving the matrices $A_1$, $A_2$, $\dots$, $A_N$.

\bd
Let $A\in M_{m,n}$.
For any index sets $\alpha$, $\beta$,
with $\alpha \subseteq\{1,\dots,m\}$,
$\beta \subseteq\{1,\dots,n\}$,
and  $|\alpha|=|\beta|$, we denote the submatrix
that lies in the rows of $A$ indexed by $\alpha$ and
the columns indexed by $\beta$ as $A(\alpha,\beta)$.
\ed
For example
$$
\left[\ba{ccccc} 1 & 0 & 9 & 0 & -2\\ 2 & 1 & 7 & 1 & 1
\ea \right](\{1\},\{1,3\})=
\left[1 \ \ 9\right].
$$

On the other hand and taking into account some
properties of the determinant it is well-known that
the characteristic polynomial of a given square matrix
$A$ can be written as
$$
\chi_A(t)=\det(t I-A)=t^n-S_1(A)t^{n-1}+\cdots
+(-1)^nS_n(A),
$$
where $I\in M_n$ is the identity, and $S_k(A)$ is the
elementary symmetric function associated to the matrix
$A$, $k=1,2,\dots,n$.

In fact, by the second way as we have defined the
determinant, i.e. the alternating sum, it is
straightforward that
\bq \label{ex4}
S_k(A)=\sum_{|\alpha|=k}\det(A(\alpha,\alpha)),\qquad
k=1,2, \dots,n.
\eq
In connection with the elementary symmetric functions
we pre\-sent new equalities related to these functions,
giving explicit expressions for $S_2(A+B)$ and
$S_3(A+B)$, for any $A, B\in M_n$.

The structure of this paper is the following: In Section
2 we  present some results related with the determinant
of sum of matrices, whose proof is given in the appendix.
In Section 3 obtain the values of $S_2(A+B)$ and $S_3(A+B)$
by using the definition of the elementary symmetric functions
of a matrix, in Section 4 we prove the same identities and
also we obtain $S_4(A+B)$ by using the Newton-Girard identities,
where $A$ and $B$ are two generic $n$-by-$n$ matrices.

\section{The determinant of a sum of matrices}

Let $N$ be a positive integer and let us
consider the $N$-tuple of $n$-by-$n$ matrices
$$
S:=(A_1, A_2, \cdots, A_N).
$$
We define $\Sigma(S)$ as the set of all possible formal
sums of matrices of $S$ where each $A_i$, $i=1,\dots, N$, appears
at most once.
\bn
Note that WLG we can add the null
matrix, $0$, in $\Sigma(S)$.
\en

The following result will be useful for further results
\bt \label{the1}
Given $A\in M_n$ and an integer $N$, with $N\ge n+1$.
For any $N$-tuple $S=(A_1,A_2,\cdots,A_N)$,
$A_i\in M_n$, $i=1,\dots,N$, the following relation holds:
\bq \label{ex1}
\sum_{k=0}^N (-1)^k \sum_{\substack{\Omega\in \Sigma(S)
\\|\Omega|=k}}\det \left(A+\sum_{A_i\in \Omega}A_i\right)=0,
\eq
understanding that $|\Omega|=k$ means that $\Omega$ is
a formal sum with $k$ summands, and that $A_i\in \Omega$
means that $A_i$ is a summand in $\Omega$.
\et

\bn
The identity \eqref{ex1} can be rewritten as
\bq \label{ex1.1}
\sum_{x_1,\dots,x_N=0}^1
(-1)^{x_1+\cdots+x_N} \det\left(
A+\sum_{j=1}^N x_j A_j\right)
=0.
\eq
Chapman proves in \cite{ams02} the case $A=0$ of this.
But his argument works as well in this generalized form;
the determinant is a polynomial of degree less that $N$
in the variables $x_1$, $\dots$, $x_N$ and this
alternating sum must vanish as seen by applying
to any monomial of degree less than $N$.
Alternatively \eqref{ex1.1} follows by subtracting the
$N+1$ case of Chapman's identity from the $N$ case.
\en
For instance, if we set  $A=0$ in \eqref{ex1} and
$A_1,A_2,A_3,A_4\in M_3$, i.e. $N=4$,
then
$$\ba{rl}
\det(A_1+A_2+A_3+A_4)=& \det(A_1+A_2+A_3)+\det(A_1+A_2+A_4)
\\ & +\det(A_1+A_3+A_4)+\det(A_2+A_3+A_4)\\ & -\det(A_1
+A_2)-\det(A_1+A_3)-\det(A_1+A_4)\\
& -\det(A_2+A_3)-\det(A_2+A_4)-
\det(A_3+A_4)\\& +\det(A_1)+\det(A_2)+\det(A_3)+
\det(A_4).
\ea $$
This result has very interesting consequences.
\bc \label{cor1}
Under the conditions of Theorem \ref{the1}.
For any index sets $\alpha$, $\beta \subseteq\{1,2,
\dots,n\}$ of size  $\tau$, $N\ge \tau+1$, the following
relation holds:
\bq \label{ex2}
\sum_{k=0}^N (-1)^k \sum_{\substack{\Omega\in \Sigma(S)
\\|\Omega|=k}}\det \left(A(\alpha,\beta)+\sum_{A_i\in \Omega}
A_i(\alpha,\beta)\right)=0.
\eq
\ec
The proof follows from Theorem \ref{the1}
replacing $A$ by $A(\alpha,\beta)$ and taking into account
that $A_i(\alpha,\beta)\in M_{\tau}$ and $N\ge \tau+1$.
\\
On the other hand, if we combine the above result
and \eqref{ex4} we obtain that:
\bc \label{cor25}
Under the conditions of Theorem \ref{the1}.
For any nonnegative integer $\tau$, $N\ge \tau+1$,
\bq \label{ex3}
\sum_{k=0}^N (-1)^k \sum_{\substack{\Omega\in \Sigma(S)
\\|\Omega|=k}}S_\tau \left(A+\sum_{A_i\in \Omega}A_i\right)=0,
\eq
where $S_\tau(C)$ is the $\tau$-th elementary symmetric
function of the matrix $C$.
\ec
The proof, again, is straightforward taking into account
\eqref{ex4} and that $N\ge \tau+1$.

The following identity is useful to compute $\tau$-th
elementary symmetric function of any number of matrices
$N\ge \tau+1$.
\bc \label{cor26}
Under the conditions of Theorem \ref{the1}.
For any nonnegative integers $\tau$, $N\ge \tau+1$, the following
identity fulfills
\bq \label{aux19}
S_\tau(A_1+A_2+\cdots+A_N)=\sum_{j=0}^{\tau-1}(-1)^{j}
\binom{j+N-\tau-1}{N-\tau-1}\sum_{\substack{\Omega\in 
\Sigma(S)\\|\Omega|=\tau-j}}S_\tau\left(\sum_{i\in 
\Omega}A_i\right).
\eq
\ec
which proof is elementary and we leave it for the reader.

In fact Theorem \ref{the1} is optimal with respect
to the range of $N$, i.e. for every positive integer $n$,
it is possible to find $n$-tuples of $M_n$ such  that
the equality \eqref{ex1}, given in Theorem \ref{the1},
fails.
For instance, taking
$$
A_i={\rm diag}(e_i),\quad i=1,2,\dots,n,\quad A=xe_1,
\ x\in \XX R,
$$
where $\{e_1,e_2,\dots,e_n\}$ is the canonical basis of
$\XX R^n$, it is straightforward to check that
$$
\sum_{k=0}^n (-1)^k \sum_{\substack{\Omega\in \Sigma(S)
\\|\Omega|=k}}\det \left(A+\sum_{A_i\in \Omega}A_i\right)=
(-1)^n(1+x-x)\ne 0.
$$

\section{Obtaining $S_2(A+B)$ and $S_3(A+B)$}
So the next logical step is to get closed
expressions for the $\tau$-th elementary symmetric functions of a
sum of $N$ matrices,  with $1\le N\le \tau$.
To do that we will use the Newton-Girard formulas
for the elementary symmetric functions (see e.g.
\cite[\S 10.12]{ser}) and the definition of such
functions \eqref{ex4}.
\bn
Note that if $A$ is $n$-by-$n$, then $\det (A)=S_n(A)$,
so it is enough to obtain those identities for the
elementary symmetric functions and then apply these to
the determinant.
\en

\bl \label{lemm32}
For any $A_1, A_2\in M_n$, we get
\begin{align}
S_2(A_1+A_2)&=S_2(A_1)+S_2(A_2)+S_1(A_1)
S_1(A_2)-S_1(A_1A_2),\label{aux7} \\ \nonumber S_3(A_1
+A_2)&= S_3(A_1)+S_3(A_2)-S_1(A_1+A_2)S_1(A_1A_2)\\
\label{aux8} & +S_1(A_1)S_2(A_2)+S_1(A_2)S_2(A_1)+
S_1(A^2_1A_2)\\ \nonumber  & +S_1(A_1A^2_2),\\
\nonumber S_3(A_1+A_2+A_3)&= S_3(A_1)+S_3(A_2)+
S_3(A_3)-S_1(A_1+A_2+A_3)\\ & \times
S_1(A_1A_2+A_1A_3+A_2A_3)+S_1(A_1)\big(S_2(A_2)
\nonumber \\  \label{aux9} & +S_2(A_3)\big)+S_1(A_2)\big(
S_2(A_1)+S_2(A_3)\big)+S_1(A_3)\\ \nonumber & \times
\big(S_2(A_1)+S_2(A_2)\big)+
S_1(A^2_1A_2)\nonumber +S_1(A_1A^2_2)\\
\nonumber & +S_1(A^2_1A_3)+S_1(A_1A^2_3)
+S_1(A^2_2A_3)+S_1(A_2A^2_3)\\
\nonumber & +S_1(A_1A_2A_3)+S_1(A_1A_3A_2).
\end{align}
\el
\bdm
Let $A_1, A_2\in  M_n$ be two matrices with
spectra $\sigma(A_1)=\{\lambda_1, \dots,
\lambda_n\}$ and $\sigma(A_2)=\{\mu_1, \dots, \mu_n\}$,
respectively.
The Newton-Girard formula gives
$$
\sum_{i=1}^n \lambda_i^2=S_1^2(A_1)-2S_2(A_1),
$$
where $S_1(A_1)\equiv Tr(A_1)$.
WLG we can assume $A_1$ diagonal, then
by definition of $S_2$ (see \eqref{ex4}) we get
$$\ba{rl}
S_2(A_1+A_2)=& \dst \sum_{1\le i<j\le n}\lambda_{i}
\lambda_{j}+2\sum_{j=1}^n \lambda_j(Tr(A_2)-b_{j,j})
+S_2(A_2)\\[4mm]
= & \dst \half \Big(\sum_{j=1}^n \lambda_j\Big)^2-
\half \sum_{j=1}^n \lambda_j^2+ 2Tr(A_1)Tr(A_2)-2Tr(
A_1A_2)\\[4mm] & +S_2(A_2)=\dst \half (Tr(A_1))^2-
\half\Big((Tr(A_1))^2-2S_2(A_1)\Big)\\[4mm]& \dst +
2Tr(A_1)Tr(A_2)-2Tr(A_1A_2)+S_2(A_2)\\[4mm] &
\dst = S_2(A_1)+S_2(A_2)+Tr(A_1)Tr(A_2)-Tr(A_1A_2).
\ea$$
And for $S_3(A_1+A_2)$, if $A_1$ is a diagonal matrix,
by definition of $S_3$ (see \eqref{ex4}) it is
straightforward to get
$$\ba{rl}
S_3(A_1+A_2)=& \hspace{-2mm}\dst \sum_{1\le i<j<k\le n}
\lambda_{i}\lambda_{j}\lambda_{k}+\sum_{1\le i<j\le n}
\lambda_{i}\lambda_{j}\Big(Tr(A_2)-b_{i,i}-b_{j,j}\Big)
\\[4mm] +& \dst \sum_{j=1}^n \lambda_j S_2((A_2)_{j,j})+
S_3(A_2).
\ea$$
Taking into account that in this case the Newton-Girard
formula produces the identity
$$
\sum_{i=1}^n \lambda_i^3=S^3_1(A_1)-3S_1(A_1)S_2(A_1)+
3S_3(A_1),
$$
and the expansion of $(a+b+c)^3$, we obtain
$$\ba{rl}
S_3(A_1+A_2)=& \dst \frac 16\Big((Tr(A_1))^3+2Tr(A_1^3)-
3Tr(A_1^2)Tr(A_1))\Big)\\+ & \dst S_2(A_1)Tr(A_2)-Tr(A_1)
Tr(A_1A_2)+Tr(A_1^2A_2)\\ +& \dst \sum_{j=1}^n \lambda_j
S_2((A_2)_{j,j})+S_3(A_2).
\ea$$
But we can assume $A_2$ is a diagonal matrix and say
$$
S_2((A_2)_{j,j})=S_2(A_2)-\mu_jTr(A_2)+\mu_j^2, \qquad
j=1, 2, \dots, n.
$$
So,
$$\ba{rl}
S_3(A_1+A_2)= & \dst S_3(A_1)+S_2(A_1)Tr(A_2)-Tr(A_1)Tr(A_1
A_2)\\+ & \dst Tr(A_1^2A_2)+3Tr(A_1)S_2(A_2)-Tr(A_1A_2)
Tr(A_2)\\+ & \dst Tr(A_1A_2^2)+S_3(A_2),
\ea$$
and hence both relations, \eqref{aux7} and \eqref{aux8}, hold.
Moreover \eqref{aux9} is a direct consequence of \eqref{aux8}.
\edm

\section{Other way to obtain $S_2(A+B)$, $S_3(A+B)$ and $S_4(A+B)$.}

We will start probing $S_2(A+B)$ using the Newton-Girard identities:
$$\ba{rl}
-2S_2(A+B)=& S_1((A+B)^2)-S_1^2(A+B)=-S_1^2(A)+S_1(A^2)-S_1^2(B)\\[3mm]
& +S_1(B^2)-2S_1(A)S_1(B)+2S_1(AB)=-2S_2(A)-2S_2(B)\\
& -2S_1(A)S_1(B)+2S_1(AB).
\ea$$
We will apply an analogous way to obtain $S_3(A+B)$:
$$\ba{rl}
3S_3(A+B)=& S_1((A+B)^3)-S_1((A+B)^2)S_1(A+B)+S_1(A+B)S_2(A\\ & +B)
=S_1^3(A)+3S_1(A^2B)+3S_1(AB^2)+S_1(B^3)-S_1(A^2)\\ & \times S_1(A)
-S_1(A^2)S_1(B)-2S_1(AB)S_1(A)-2S_1(AB)S_1(B)\\ & -S_1(B^2)S_1(A)
-2S_1(B^2)S_1(B)+S_1(A+B)S_2(A+B).
\ea$$
If now we expand $S_2(A+B)$, after some simplifications it is
clear we get the desired identity for $S_3(A+B)$.
\subsection{Obtaining $S_4(A+B)$} \label{4.1}
As the above examples, the Newton-Girad formula gives
$$
\ba{rl}
-4S_4(A+B)=& S_1\left((A+B)^4\right)-S_1\left(
(A+B)^3\right)S_1(A+B)+S_1\left((A+B)^2\right)
\\[3mm] & \times S_2(A+B)-S_1(A+B)S_3(A+B).
\ea$$
Taking into account the properties of the trace,
we get
$$
\ba{rl}
-4S_4(A+B)=& \hspace{-3mm} S_1(A^4)+4S_1(A^3B)+4
S_1(A^2B^2)+2S_1\left((AB)^2\right)+4
S_1(AB^3)\\[2mm] & \hspace{3mm} +S_1(B^4)-S_1(A^3)S_1(A)-S_1(A^3)
S_1(B)-3S_1(A^2B)S_1(A)\\[3mm] & \hspace{3mm} -3S_1(A^2B)S_1(B)
-3S_1(AB^2)S_1(A)-3S_1(AB^2)S_1(B)\\[3mm] & \hspace{3mm}
-S_1(B^3)S_1(A)-S_1(B^3)S_1(B)+\big(S_1(A^2)+2S_1(AB)
\\[3mm] & \hspace{3mm} +S_1(B^2)\big)\big(S_2(A)+S_2(B)+S_1(A)S_1(B)-
S_1(AB)\big)\\[3mm] & \hspace{3mm} -(S_1(A+B))(S_3(A)+S_3(B)
+S_1(A)S_2(B)+S_1(B)S_2(A)\\[3mm] & \hspace{3mm}  +S_1(A^2B)+
S_1(AB^2)-S_1(AB)S_1(A)-S_1(AB)S_1(B))
\ea$$
Applying the same technique applied before, we get
$$
\ba{rl}
-4S_4(A+B)=& \hspace{-2mm}
-4S_4(A)-4S_4(B)
\\ & \hspace{-2mm} 4S_1(A^3B)+4S_1(A^2B^2)+2S_1\left((AB)^2\right)+
4S_1(AB^3)\\[3mm] & \hspace{-2mm} -S_1(A^3)S_1(B)-
3S_1(A^2B)S_1(A)-3S_1(A^2B)S_1(B)\\[3mm] & \hspace{-2mm} -3S_1(AB^2)S_1(A)-
3S_1(AB^2)S_1(B)-S_1(B^3)S_1(A)\\[3mm] & \hspace{-2mm}
+S_1(A^2)(S_2(B)+S_1(A)S_1(B)-S_1(AB))
+2S_1(AB)(S_2(A)\\[3mm] & \hspace{-2mm} +S_2(B)+S_1(A)S_1(B)-S_1(AB))
+S_1(B^2)(S_2(A)\\[3mm] & \hspace{-2mm} +S_1(A)S_1(B)-S_1(AB))
-S_1(A)(S_3(B)+S_1(A)S_2(B)\\[3mm] & \hspace{-2mm}+S_1(B)S_2(A)+S_1(A^2B)+
S_1(AB^2)\\[3mm] & \hspace{-2mm} -S_1(AB)S_1(A)-S_1(AB)S_1(B))-S_1(B)(S_3(A)
\\[3mm] & \hspace{-2mm}+S_1(A)S_2(B)+S_1(B)S_2(A)+S_1(A^2B)+S_1(AB^2)\\[3mm]
& \hspace{-2mm} -S_1(AB)S_1(A)-S_1(AB)S_1(B)).
\ea$$
After some simplifications applying the Newton-Girard formulas
we get
$$
\ba{rl}
4S_4(A+B)=&\hspace{-2mm} 4S_4(A)+4S_4(B)-4S_1(A^3B)-4S_1(A^2B^2)-
2S_1\left((AB)^2\right)\\[3mm] &\hspace{-2mm}-4S_1(AB^3) +3S_1(A^2B)
S_1(A)+3S_1(A^2B)S_1(B) +3S_1(AB^2)\\[3mm] &\hspace{-2mm} \times S_1(A)+
3S_1(AB^2)S_1(B)-2S_1(AB)(S_1(A)S_1(B)\\[3mm] &\hspace{-2mm}-S_1(AB))
+S_1(A)(S_1(A^2B)+S_1(AB^2)-S_1(AB)
S_1(B))\\[3mm] &\hspace{-2mm}+S_1(B)(S_1(A^2B)+S_1(AB^2)-S_1(AB)
S_1(A))\\[3mm] &\hspace{-2mm} +4S_3(A)S_1(B) +4S_3(B)S_1(A)+4S_2(A)
S_2(B)\\[3mm] &\hspace{-2mm} -4S_2(A)S_1(AB)-4S_2(B)S_1(AB).
\ea$$
Applying the Newton-Girard formulas and
after some simplifications, we get
$$
\ba{rl}S_4(A+B)=&\hspace{-2mm} S_4(A)+S_4(B)-S_1(A^3B)-S_1(A^2B^2)-
S_1(AB^3)+S_1(A^2B)\\[3mm] &\hspace{-2mm} \times S_1(A)+S_1(A^2B)S_1(B)+S_1(AB^2)
S_1(A)+S_1(AB^2)S_1(B)\\[3mm] &\hspace{-2mm} -S_1(AB)S_1(A)S_1(B)
+S_3(A)S_1(B)+S_3(B)S_1(A)+S_2(A)\\[3mm] &\hspace{-2mm} \times S_2(B)
-S_2(A)S_1(AB)-S_2(B)S_1(AB)+S_2(AB).
\ea$$

\section{Conclusions and Outlook}

We have constructed the 2nd, the 3rd and the 4th elementary
symmetric function of a sum of two matrices but, of course,
is simple to see that is possible to compute the $\tau$-th
elementary symmetric function of a sum of $N$-matrices, $1\le N\le
\tau$ by using the Newton-Girard formulas or by using the same
technique used in Lemma \ref{lemm32} which, by the way, is too
much complicated.

Of course, one of the goals in further papers is to find a
closed expression in the general case which for the moment
is not clear although we believe the Theory of partition of
integers is involved.

In fact, by using the generalized Waring's formula \cite{zhe},
for any $A$ is $n$-by-$n$ matrix, $0\le k\le n$ and any
nonnegative integer $m$, we get
$$
S_k(A^m)=(-1)^{k(n+1)}\sum_{|\lambda|=kn} A_\lambda e_\lambda(A),
$$
where the coefficients $A_{\lambda}$ are given by
$$
A_\lambda=\sum_{\substack{|\pi|=k \\ \pi=(k_1,k_2,\dots)}}
\sum_{\substack{\lambda_1\cup \lambda_2\cup \cdots \cup \lambda_{l(\pi)=\lambda}
\\|\lambda_i|=k_in}}\frac{(-1)^{l(\lambda)-l(\pi)}}{\prod_{i=1}^k m_i(\pi)!}
\prod_{l=1}^{l(\pi)}\frac n {l(\lambda_i)}\binom{l(\lambda_i)}{m_1(\lambda_i),\dots,
m_m(\lambda_i)},
$$
and  $e_\lambda(A)=S_1^{m_1(\lambda)}(A)
S_2^{m_2(\lambda)}(A)\cdots$

\bn
A {\textit{partition}}is a finite sequence $(\lambda_1,\lambda_2,\dots,
\lambda_r)$ of positive integers in decreasing order, where
$l(\lambda)$ denotes the length of the partition, and $m_k(\lambda)$ denotes
the number of parts of $\lambda$ equal to $k$.

$\lambda\cup \mu$ is the partition whose parts are those of $\lambda$ and $\mu$.
\en

Taking into account this identity we believe that we can obtain an
analogous expression for the $\tau$-th elementary symmetric function
of a sum of matrices. In fact, we expect one expression
in which appears the elementary symmetric functions on words of the
letters $A_1$, $A_2$, $\dots$, $A_N$ as one could see in  subsection
\ref{4.1} for the case $S_4(A_1+A_2)$.

{\bf Acknowledgements:}
The author thanks the referee for the constructive
remarks and the valuable comments.
This work has been supported by Direcci´on General de
Investigaci´on (Ministerio de Educaci´on y Ciencia) of
Spain, grant MTM 2006-13000-C03-02.

\appendix
\section{Proof of Theorem \ref{the1}}
We will prove by induction on $n$:
\bi
\item If $n=1$ the matrices are scalars so, for every
$k$,
$$\ba{rl}
\dst \sum_{\substack{\Omega\in \Sigma(S)\\|\Omega|=k}}\det
(A+&\sum_{A_i\in \Omega}A_i)= \dst \sum_{\substack{\Omega
\in \Sigma(S) \\|\Omega|=k}} (A+\sum_{A_i\in \Omega}A_i)
\\[8mm] =&  \dst \binom N k A+\binom {N-1}{k-1}(A_1+A_2
+\cdots+A_N),
\ea$$
and hence \eqref{ex1} holds for $n=1$ and $N\ge 2$.
\item If we assume that the result holds for $n$, let us
going to prove  the identity \eqref{ex1} for $n+1$:
\\
Taking the Laplace expansion through the first row, we get
$$
\ba{c}
\dst\sum_{k=0}^N (-1)^k \sum_{\substack{\Omega\in\Sigma(S)
\\|\Omega|=k}}\det (A+\sum_{A_i\in \Omega}A_i)=\sum_{k=0}^N
(-1)^k \sum_{\substack{\Omega\in
\Sigma(S) \\|\Omega|=k}}\Big\{\sum_{j=1}^{n+1} (-1)^{1+j}
\\[5mm] \dst \times \Big(A(1,j)+\sum_{A_i\in \Omega}A_i(1,
j)\Big)\det\Big(A_{1,j}+\sum_{A_i\in \Omega} (A_i)_{1,
j}\Big)\Big\}.
\ea
$$
By induction, since for every $j=1, 2, \dots, n$,
$A(1,j)$ is fixed and does not depend on $k$ nor $\Omega$
we get that the above expression is equal to
$$
\ba{c} \dst \sum_{j=1}^{n+1} (-1)^{j+1}\sum_{k=0}^N
(-1)^k \sum_{\substack{\Omega\in \Sigma(S)\\|\Omega|=k}}
\Big(\sum_{A_i\in \Omega}A_i(1,j)\Big)\det\Big(A_{1,j}+
\sum_{A_i\in \Omega}(A_i)_{1,j}\Big)\\[5mm]
\dst =\sum_{j=1}^{n+1}(-1)^{j+1}\sum_{k=0}^N (-1)^k
\sum_{\lambda=1}^N A_\lambda(1,j) \sum_{\substack{\Omega\in
\Sigma(S) \\A_\lambda \in\Omega \\|\Omega|=k}}\det(A_{1,j}
+ \dst \sum_{A_i\in \Omega}(A_i)_{1,j})
\ea
$$
Now, if we assume that any set with less than one element
has determinant equal to zero, we get
$$
\ba{c} \dst \sum_{j=1}^{n+1}(-1)^{j+1} \sum_{k=0}^N
(-1)^k\sum_{\lambda=1}^N A_\lambda(1,j) \hspace{-3mm}
\sum_{\substack{\widetilde \Omega\in \Sigma(S\setminus\{
A_\lambda\})\\|\widetilde \Omega|=k-1}}\hspace{-3mm}\det(
A_{1,j}+(A_\lambda)_{1,j}+\sum_{A_i\in \widetilde \Omega}
(A_i)_{1,j}) \\[5mm] = \dst \sum_{j=1}^{n+1}
(-1)^{j+1} \sum_{\lambda=1}^N A_\lambda(1,j)\sum_{k=0}^{N-
1} (-1)^k \hspace{-4mm}\sum_{\substack{\widetilde
\Omega\in \Sigma(S\setminus\{A_\lambda\})\\|\widetilde
\Omega|=k}} \hspace{-4mm}\det(A_{1,j}+ (A_\lambda)_{1,j}+
\sum_{A_i\in \widetilde \Omega}(A_i)_{1,j}) \\[5mm]
\ea
$$
By induction, since for every $\lambda$ and $j$,
the matrices $A_{1,j}, (A_\lambda)_{1,j}\in M_{n-1}$ are
fixed, thus $A=A_{1,j}+(A_\lambda)_{1,j}\in M_{n-1}$ is
also fixed.
Thus we get
$$
\sum_{j=1}^{n+1}(-1)^{j+1}\sum_{\lambda=1}^N
A_\lambda(1,j) \cdot 0=0.
$$
Moreover, since $N\ge n+2$, then $N-1\ge n+1$.
Hence, the relation holds.

\hfill \rule{2ex}{2ex}
\ei


\begin{thebibliography}{99}
\bibitem{aig79} M. Aigner, {\it Combinatorial Theory},
Springer-Verlag, originally published 1979 as Vol. 234
of the Grundlehren der math. Wissenchaften, reprinted 1997.

\bibitem{BeMoVi75} D. Bessis, P. Moussa and M. Villani:
Monotonic converging variational approximations to
the functional integrals in quantum statistical mechanics,
J. Math. Phys. {\bf 16} (1975), 2318–2325.
\bibitem{HoJo85} R. Horn and C. Johnson, {\it Matrix
analysis}, Cambridge University Press (Cambridge, 1985).

\bibitem{LiSe04} E.H. Lieb and R. Seiringer: Equivalent
forms of the Bessis-Moussa-Villani conjecture, J.
Statist. Phys. {\bf 115} (1--2) (2004), 185–-190.

\bibitem{ser} S\'{e}roul, R., {\it Programming for
Mathematicians}, Springer-Verlag (Berlin, 2000).

\bibitem{zhe} Jiang Zheng: On a generalization of Waring's formula,
Advanc. Appl. math. {\bf 19} (1997), 450--452.

\bibitem{ams02} Amer. Math. Monthly {\bf 109}(7)
(2002), 665--666.
\end{thebibliography}
\end{document}